\documentclass[12pt]{article}

\usepackage{epsfig}
\usepackage{amssymb}
\usepackage{amsmath}
\usepackage{amsfonts}
\usepackage{bbm}
\usepackage{hyperref}

\def\diam{{\rm diam}}

\newtheorem{theorem}{Theorem}[section]

\newtheorem{proposition}[theorem]{Proposition}
\newtheorem{definition}[theorem]{Definition}

\newtheorem{remark}[theorem]{Remark}

\title{Off-Criticality and the \\ Massive Brownian Loop Soup}

\author{
{Federico Camia}
\thanks{NYU Abu Dhabi and VU University Amsterdam. E-mail: federico.camia\,@\,nyu.edu}
}

\begin{document}

\date{}

\maketitle

\begin{abstract}
We introduce a natural ``massive'' version of the Brownian loop soup of Lawler and Werner which
displays conformal covariance and exponential decay. 
We show that this massive Brownian loop soup arises as the near-critical scaling limit of a random
walk loop soup with killing and is related to the massive SLE$_2$ identified by Makarov and Smirnov
as the near-critical scaling limit of a loop-erased random walk with killing. We conjecture that the
massive Brownian loop soup describes the zero level lines of the massive Gaussian free field, and
discuss possible relations to other models, such as Ising, in the off-critical regime.
\end{abstract}

\medskip\noindent
{\em Key words and phrases}: Brownian loop soup, random walk loop soup, off-critical regime,
near-critical scaling limit, conformal covariance.

\medskip\noindent
{\em AMS subject classification}: 60J65, 82B41, 60K35, 81T40.



\section{Introduction and Motivation} \label{overview}
The Brownian loop soup can be thought of, roughly speaking, as an ``ideal gas'' of planar Brownian loops,
or an ensemble of loops obtained as a Poisson realization of a conformally-invariant measure on loops
constructed starting from the Brownian bridge measure.
It was introduced and studied by Lawler and Werner \cite{lw} because of its conformal invariance
properties and its close connections to the Schramm-Loewner Evolution (SLE) and to some Conformal Loop
Ensembles (CLEs) \cite{sw}. Among other properties, in some cases it can be used to describe the Radon-Nikodym
derivative between the measures corresponding to (chordal) SLE curves in different domains
(see, e.g., \cite{lawler09} and \cite{lawler-survey}).

The goal of this paper is to introduce a natural ``massive'' extension of the Brownian loop soup that doesn't
seem to have been explicitly defined and studied before, except for the lecture notes \cite{camia-notes},
which contain many of the results presented here, but in a more diluted form and with a different emphasis.
Here we focus on the relation of the massive Brownian loop soup, as we call it, to the off-critical regime of
models such as the random walk loop soup \cite{ltf}, the loop-erased random walk and the Gaussian free field.
The results we present in this paper are obtained using standard techniques and are relatively simple to prove.
Nonetheless, we believe that they are interesting because they suggest that this extension of the Brownian
loop soup may play a significant role in the study of near-critical scaling limits, much as the ``critical''
(i.e., conformally-invariant) Brownian loop soup introduced in \cite{lw} plays a role in the description of
critical scaling limits.

The statement above is motivated in part by the conformal covariance property of the massive Brownian loop soup,
expressed by Equation \eqref{eq:conf-cov} in Section~\ref{sec:bls}, by the exponential decay of its loop clusters
(second bullet of Theorem~\ref{thm:phase-trans}), and by the fact that the massive Brownian loop soup arises as
the near-critical scaling limit of the random walk loop soup (see Section \ref{sec:rwls} and particularly Theorem
\ref{prop:massive-soups}). But maybe more convincing is Theorem~\ref{theorem-rd-derivative} discussed below,
which links the massive Brownian loop soup to a massive version of SLE$_2$ identified by Makarov and Smirnov \cite{ms10},
as we now explain.

Suppose that $D$ is a bounded, simply connected domain with 
two marked points, $a,b \in \partial D$, on its boundary. Lawler, Schramm and Werner showed \cite{lsw04}
that the loop-erased random walk from $a$ to $b$ in $D$ converges in the scaling limit to chordal SLE$_2$
in $D$ from $a$ to $b$.

Being ``built from'' the simple symmetric random walk and having a conformally-invariant scaling limit,
the loop-erased random walk can be seen as a critical model. A natural non-critical version of the loop-erased
random walk is obtained by replacing the simple symmetric random walk by a random walk with killing. 
In \cite{ms10}, Makarov and Smirnov identified various massive versions of SLE, one of which, the massive SLE$_2$,
corresponds to the near-critical scaling limit of the loop-erased random walk with killing (see Section 2.1 of \cite{ms10}).

Now suppose that $D' \subset D$ is a second simply connected
domain 
such that $\partial D$ and $\partial D'$ agree in a neighborhood of both $a$ and $b$. Let ${\cal P}_{D}^{m}$
denote the distribution of the massive SLE$_2$ in $D$ from $a$ to $b$, and ${\cal P}_{D'}^{m}$ the distribution
of the massive SLE$_2$ in $D'$ from $a$ to $b$ (for ease of notation, we omit the dependence on $a,b$, which we assume
fixed). Let also ${\mathbb P}^{D}_{\lambda,m}$ denote the distribution of a massive Brownian loop soup in $D$ with
intensity $\lambda$. 
(Like the ``critical'' Brownian loop soup, the massive Brownian loop soup is characterized by an intensity $\lambda>0$.)
Then, assuming Theorem~2.1 of \cite{ms10}, one has the following result (where $\mathbbm{1}_{\{\cdot\}}$
denotes the indicator function).
\begin{theorem} \label{theorem-rd-derivative}
${\cal P}_{D'}^{m}$ is absolutely continuous with respect to ${\cal P}_{D}^{m}$
and the Radon-Nikodym derivative is
\begin{equation} \label{intro-rd-derivative}
\frac{d{\cal P}_{D'}^{m}}{d{\cal P}_{D}^{m}}(\eta) = \frac{Y_{m}(\eta)}{E_{D}^{m}(Y_{m})} \, ,
\end{equation}
where
\begin{equation} \nonumber
Y_m (\eta) = Y_m(D,D';a,b)(\eta) := {\mathbb P}^{D}_{1,m}(\text{no loop intersects } D \setminus D' \text{ and } \eta)
\mathbbm{1}_{\{ \eta \subset D' \}} ,
\end{equation}
and $E_{D}^{m}$ denotes expectation with respect to ${\cal P}_{D}^{m}$.
\end{theorem}

The analog of Eq.~\eqref{intro-rd-derivative} involving the corresponding massless objects is well known
(see, e.g., Section~6 of \cite{lawler09}, and Section~6.7 of~\cite{lawler-survey}). The fact that such a relation holds in the
massive case motivates us to study the massive Brownian loop soup with distribution ${\mathbb P}^{D}_{\lambda,m}$ that
appears in Theorem~\ref{theorem-rd-derivative}. 

In the massless case, as seen in \cite{lawler09} and \cite{lawler-survey}, the relation expressed by Eq.~\eqref{intro-rd-derivative}
holds for general SLE$_{\kappa}$ curves with $0 < \kappa \leq 8/3$, with $\lambda=1$ replaced by
$\lambda=\lambda(\kappa)=-\frac{(3\kappa-8)(6-\kappa)}{4\kappa}$.\footnote{We note that there is some
confusion in most of the existing literature regarding the relation between $\lambda$ and $\kappa$. This is due
to the fact that the Brownian loop measure used as intensity measure in defining the Brownian loop soup is an
infinite measure. As a consequence, a choice of normalization is required when defining the Brownian loop soup.
This choice is then reflected in the relation between $\lambda$ and $\kappa$. Our choice of normalization is
made explicit in Section \ref{sec:bls} and coincides with that of \cite{lw}.
We thank Greg Lawler for a useful discussion on this topic.}
It is tempting to conjecture that a similar relation holds in such generality also in the massive case. One problem with
making such a conjecture precise is that, while the massive Brownian loop soup of intensity $\lambda$ is a uniquely
defined object, in general there can be more than one way to perturb an SLE curve to obtain a massive version.

If we introduce $c(\kappa) = \frac{(3\kappa-8)(6-\kappa)}{2\kappa}$, in the massless case, $m \equiv 0$, dropping the
$m$ from the notation, we can write Eq.~\eqref{intro-rd-derivative} as
\begin{eqnarray*}
\lefteqn{
\frac{d{\cal P}_{D'}}{d{\cal P}_{D}} (\eta) } \\
& = & \frac{1}{E_{D}(Y)} \, \exp\left(\frac{c(\kappa)}{2}\mu_D(A(D;D',\eta))\right) \, \mathbbm{1}_{\{ \eta \subset D' \}} ,
\end{eqnarray*}
where $A(D;D',\eta) := \{\text{loops in } D \text{ that intersect } D \setminus D' \text{ and } \eta \}$,
$\mu_D$ is the intensity measure of the Brownian loop soup (the Brownian loop measure defined in the next section
and studied in \cite{werner3}), and $Y = Y(D,D';a,b)(\eta) := \exp\left(\frac{c(\kappa)}{2}\mu_D(A(D;D',\eta))\right)
\mathbbm{1}_{\{ \eta \subset D' \}}$.
This expression is valid for all values of $\kappa$ between $0$ and $4$ (see, e.g., Section~6 of \cite{lawler09} and
Section~6.7 of~\cite{lawler-survey} again).
It is again tempting to think that such a relation may hold in the massive case, with $\mu_{D}$ replaced by the massive Brownian
loop measure $\mu^{m}_{D}$ introduced in the next section. Once again though, to make this conjecture precise one would need
to identify a specific massive version of SLE.


Based on the considerations above and on \cite{bcl14,ww14,lupu15}, a precise conjecture can be made in the context of the
Gaussian free field.\\

\noindent{\bf Conjecture.} The collection of zero level lines of the massive Gaussian free field in a bounded, planar domain $D$ with
zero boundary condition is distributed like the collection of cluster boundaries from a massive Brownian loop soup in $D$ with intensity
$\lambda=1/2$.\\

We remark that it is plausible that an extension of the methods of Lupu \cite{lupu15} will lead to a proof of the conjecture\footnote{The
author thanks Titus Lupu for a useful discussion concerning \cite{lupu15}}.

One could attempt to formulate a similar conjecture in the context of the Ising model by saying that the collection of spin-cluster boundaries
in the near-critical scaling limit is distributed like the collection of cluster boundaries from a massive Brownian loop soup with intensity $\lambda=1/4$.
The obvious problem with that formulation is that there are at least two distinct, natural ways to take a near-critical scaling limit in the context
of the Ising model: by adding an external magnetic field at the critical temperature (see \cite{cgn13}), and by moving away from the critical
temperature at zero external magnetic field. These two near-critical scaling limits are believed to lead to different massive versions of SLE$_3$,
highlighting the problem discussed above. We believe that the massive Brownian loop soup with intensity $\lambda=1/4$ is related to one of
the two near-critical scaling limits described above, but at the moment we don't have any convincing evidence for picking one of the two and
are therefore not able to formulate a precise conjecture for the Ising model.\\

%

\noindent\emph{Organization of the Paper.} In Section \ref{sec:bls}, we introduce the massive Brownian loop soup and present some
of its properties. In Section \ref{sec:rwls}, we introduce the random walk loop soup and study its near-critical scaling limit, showing
that it leads to the massive Brownian loop soup. Finally, in Section \ref{sec:boundary}, we prove Theorem~\ref{theorem-rd-derivative}.

\section{Brownian Loop Soups} \label{sec:bls}


The Brownian loop soup was introduced by Lawler and Werner in \cite{lw}. Following \cite{lw},
we call \emph{rooted loop} a continuous function $\gamma:[0,t_{\gamma}] \to {\mathbb C}$ with
$\gamma(0)=\gamma(t_{\gamma})$. We will consider only loops with $t_{\gamma} \in (0,\infty)$.
The \emph{Brownian bridge measure} $\mu^{br}$ is the probability measure on rooted loops of duration 1
with $\gamma(0)=0$ induced by the Brownian bridge $B_t := W_t - t W_1$, $t \in [0,1]$, where $W_t$
is standard, two-dimensional Brownian motion. A measure $\mu^{br}_{z,t}$ on loops rooted at
$z \in {\mathbb C}$ (i.e., with $\gamma(0)=z$) of duration $t$ is obtained from $\mu^{br}$ by
Brownian scaling, using the map
\begin{equation} \nonumber
(\gamma,z,t) \mapsto z + t^{1/2} \gamma(s/t) , \; s \in [0,t] \, .
\end{equation}
More precisely, we let
\begin{equation} \label{eq:brownian-bridge}
\mu^{br}_{z,t}(\cdot) := \mu^{br}(\Phi^{-1}_{z,t}(\cdot)) \, ,
\end{equation}
where
\begin{equation} \label{eq:map}
\Phi_{z,t}: \gamma(s), s \in [0,1] \mapsto z + t^{1/2} \gamma(s/t) , s \in [0,t] \, .
\end{equation}

The \emph{rooted Brownian loop measure} is defined as
\begin{equation} \label{eq:rooted-brownian-loop-measure}
\mu_r := \int_{\mathbb C} \int_0^{\infty} \frac{1}{2 \pi t^2} \, \mu^{br}_{z,t} \, dt \, d{\bf A}(z) \, ,
\end{equation}
where $\bf A$ denotes area.

The (\emph{unrooted}) \emph{Brownian loop measure} $\mu$ is obtained from the rooted one by
``forgetting the root.'' More precisely, if $\gamma$ is a rooted loop,
$\theta_u \gamma: t \mapsto \gamma(u+t \mod t_{\gamma})$ is again a rooted loop. This defines
an equivalence relation between rooted loops, whose equivalence classes we refer to as (\emph{unrooted})
\emph{loops}; $\mu(\gamma)$ is the $\mu_r$-measure of the equivalence class $\gamma$.
With a slight abuse of notation, in the rest of the paper we will use $\gamma$ to denote an unrooted loop
and $\gamma(\cdot)$ to denote any representative of the equivalence class $\gamma$.

%

The \emph{massive} (\emph{unrooted}) \emph{Brownian loop measure} $\mu^m$ is defined
by the relation
\begin{equation} \label{eq:massive-brrownian-loop-measure}
d\mu^m(\gamma) = \exp{(-R_m(\gamma))} \, d\mu(\gamma) \, ,
\end{equation}
where $m:{\mathbb C} \to {\mathbb R}$ is a nonnegative \emph{mass function} and
\begin{equation} \nonumber
R_m(\gamma) := \int_0^{t_{\gamma}} m^2(\gamma(t)) dt
\end{equation}
for any rooted loop $\gamma(t)$ in the equivalence class of the unrooted loop $\gamma$.
(Analogously, one can also define a massive \emph{rooted} Brownian loop measure:
$d\mu^m_r(\gamma) := \exp{(-R_m(\gamma))} \, d\mu_r(\gamma)$.)

If $D$ is a subset of $\mathbb C$, we let $\mu_D$ (respectively, $\mu^{m}_D$) denote $\mu$
(resp., $\mu^{m}$) restricted to loops that lie in $D$. The family of measures $\{ \mu_D \}_D$
(resp., $\{ \mu^{m}_D \}_D$), indexed by $D \subset {\mathbb C}$, satisfies the
\emph{restriction property}, i.e., if $D' \subset D$, then $\mu_{D'}$ (resp., $\mu^{m}_{D'}$) is
$\mu_D$ (resp., $\mu^{m}_D$) restricted to loops lying in $D'$.

An equivalent characterization of the Brownian loop measure $\mu$ is as follows (see~\cite{werner3}).
Given a conformal map $f:D \to D'$, let $f\circ\gamma(s)$ denote the loop $f(\gamma(t))$ in $D'$
with parametrization
\begin{equation} \label{eq:time-parametrization}
s=s(t) = \int_0^t |f'(\gamma(u))|^2 du \, .
\end{equation}
Given a subset $A$ of the space of loops in $D$, let
$f \circ A = \{ \hat\gamma = f \circ \gamma \text{ with } \gamma \in A \}$.
Up to a multiplicative constant,
$\mu_D$ is the unique measure satisfying the following two properties,
collectively known as \emph{conformal restriction property}.
\begin{itemize}
\item 
For any conformal map $f:D \to D'$,
\begin{equation} \label{eq:conf-inv}
\mu_{D'}(f \circ A) = \mu_D(A) \, .
\end{equation}
\item If $D' \subset D$, $\mu_{D'}$ is $\mu_D$ restricted to loops that stay in $D'$.
\end{itemize}

As a consequence of the conformal invariance of $\{ \mu_D \}_D$, the family of massive
measures $\{ \mu^{m}_D \}_D$ satisfies a property called \emph{conformal covariance},
defined below.
Given a conformal map $f:D \to D'$, let $\tilde m$ be defined by the map
\begin{equation} \label{eq:mass}
m(z) \stackrel{f}{\mapsto} \tilde m(w) = \left| f'(f^{-1}(w)) \right|^{-1} m(f^{-1}(w)) \, ,
\end{equation}
where $w = f(z)$. This definition, combined with \eqref{eq:time-parametrization} and
Brownian scaling, implies that $m^2 dt = \tilde m^2 ds$. From this and \eqref{eq:conf-inv},
it follows that
\begin{equation} \label{eq:conf-cov}
\mu^{\tilde m}_{D'}(f \circ A) = \mu^m_D(A) \, ,
\end{equation}
where $A$ and $f \circ A$ have the same meaning as in Eq.~(\ref{eq:conf-inv}).
We call this property conformal covariance and say that the massive Brownian
loop measure $\mu^m$ is \emph{conformally covariant}.

\begin{definition} \label{def:Bls}
A \emph{Brownian loop soup} 
in $D$ with intensity $\lambda$ is a Poissonian realization from $\lambda\mu_D$. 
A \emph{massive Brownian loop soup} 
in $D$ with intensity $\lambda$ and mass function $m$ is a Poissonian realization from $\lambda\mu^{m}_D$.
\end{definition}

The (massive) Brownian loop soup ``inherits'' the property of conformal invariance (covariance) from the
(massive) Brownian loop measure. Note that in a homogeneous massive Brownian loop soup, that is, if
$m$ is constant, loops are exponentially suppressed at a rate proportional to their time duration. We will
sometimes call the conformally invariant Brownian loop soup introduced by Lawler and Werner \emph{critical},
to distinguish it from the \emph{massive} Brownian loop soup defined above.

The definition of the massive Brownian loop soup has a nice interpretation in terms of ``killed'' Brownian motion.
For a given function $f$ on the space of loops, one can write
\begin{equation} \nonumber
\int f(\gamma) e^{-R_m(\gamma)} d\mu(\gamma) 
= \int {\mathbb E}_{T_{\gamma}} f(\gamma) \mathbbm{1}_{\{R_m(\gamma)<T_{\gamma}\}} d\mu(\gamma) \, , 
\end{equation}
where ${\mathbb E}_{T_{\gamma}}$ denotes expectation with respect to the law of the mean-one,
exponential random variable $T_{\gamma}$. 
In view of this, one can think of the Brownian loop $\gamma$ under the measure $\mu^m$
as being ``killed'' at rate $m^2(\gamma(t))$. More precisely, one has the following alternative and
useful characterization.
\begin{proposition} \label{coupling1}
A \emph{massive} Brownian loop soup in $D$ with intensity $\lambda$ and mass function $m$ can
be realized in the following way.
\begin{enumerate}
\item Take a realization of the \emph{critical} Brownian loop soup in $D$ with intensity $\lambda$.
\item Assign to each loop $\gamma$ of duration $t_{\gamma}$ an independent, mean-one,
exponential random variable, $T_{\gamma}$.
\item 
Remove from the soup all loops $\gamma$ such that
\begin{equation} \label{eq:loop-killing-1}
\int_0^{t_{\gamma}} m^2(\gamma(t)) dt > T_{\gamma} \, .
\end{equation}
\end{enumerate}

\end{proposition}
\begin{remark}
Note that Eq.~\eqref{eq:loop-killing-1} requires choosing a time parametrization
for the loop $\gamma$ but is independent of the choice.
\end{remark}

\noindent{\bf Proof.}
Let ${\cal L}^D$ denote the set of loops contained in $D$ and define
${\cal L}^D_{>\varepsilon} := \{\gamma \in {\cal L}^D : \text{diam}(\gamma)>\varepsilon \}$ and
${\cal L}^D_{>\varepsilon,r} := \{\gamma \in {\cal L}^D_{>\varepsilon} : R_m(\gamma)=r \}$.
For a subset $A$ of ${\cal L}^D_{>\varepsilon}$, let $A_r = A \cap {\cal L}^D_{>\varepsilon,r}$.
For every $\varepsilon>0$, the restriction to loops of diameter larger
than $\varepsilon$ of the massive Brownian loop soup in $D$ with mass function $m$ is a Poisson point process
on ${\cal L}^D_{>\varepsilon}$ such that the expected number of loops in $A \subset {\cal L}^D_{>\varepsilon}$
at level $\lambda>0$ is
\begin{eqnarray*}
\lambda \mu^m(A) & = & \lambda \int_A e^{-R_m(\gamma)} d\mu(\gamma) \\
& = & \lambda \int_0^{\infty} \int_{A_r} e^{-r} d\mu(\gamma) dr \\
& = & \lambda \int_0^{\infty} e^{-r} \mu(A_r) dr \, .
\end{eqnarray*}

We will now show that, when attention is restricted to loops of diameter larger than $\varepsilon$, the
construction of Proposition~\ref{coupling1} produces a Poisson point process on ${\cal L}^D_{>\varepsilon}$
with the same expected number of loops at level $\lambda>0$.

Let $N_{\lambda}(A)$ denote the number of loops in $A$ obtained from the construction of Proposition~\ref{coupling1}.
Because the Brownian loop soup is a Poisson point process and loops are removed independently, for every
$A \subset {\cal L}^D_{>\varepsilon}$ we have that
\begin{enumerate}
\item[(i)] $N_{0}(A)=0$,
\item[(ii)] $\forall \lambda,\delta>0$ and $0 \leq \ell \leq \lambda$, $N_{\lambda+\delta}(A)-N_{\lambda}(A)$ is independent of $N_{\ell}(A)$,
\item[(iii)] $\forall \lambda,\delta>0$, $\Pr(N_{\lambda+\delta}(A)-N_{\lambda}(A) \geq 2) = o(\delta)$,
\item[(iiii)] $\forall \lambda,\delta>0$, $\Pr(N_{\lambda+\delta}(A)-N_{\lambda}(A)=1) = \mu^m(A) \delta + o(\delta)$,
\end{enumerate}
where (iiii) follows from the fact that, conditioned on the event that a single additional loop appears going
from $\lambda$ to $\lambda+\delta$, the additional loops is distributed according to the density
$\frac{\mu(A_r) dr}{\mu(A)}$ on $A$.
Conditions (i)-(iiii) ensure that the point process is Poisson.

In order to identify the Poisson point process generated by the construction of Proposition~\ref{coupling1}
with the massive Brownian loop soup, it remains to compute the expected number of loops in $A$ at
level $\lambda$. For every $\varepsilon>0$ and $A \subset {\cal L}^D_{>\varepsilon}$, this is given by
\begin{equation} \nonumber
\int_0^{\infty} e^{-r} \lambda\mu(A_r) dr = \lambda\mu^m(A),
\end{equation}
which concludes the proof. \fbox{} \\

We conclude this section with a result on the connectivity phase transition and exponential decay of the
massive Brownian loop soup. Let ${\cal A}(\lambda,m,D)$ denote a massive Brownian loop soup in
$D \subseteq {\mathbb C}$ with mass function $m$ and intensity $\lambda$.
We say that two loops are \emph{adjacent} if they intersect; this adjacency relation defines \emph{clusters}
of loops, denoted by $\cal C$. (Note that clusters can be nested.) For each cluster $\cal C$, we write
$\overline{\cal C}$ for the closure of the union of all the loops in $\cal C$; furthermore, we write $\hat{\cal C}$
for the \emph{filling} of $\cal C$, i.e., the complement of the unbounded connected component of
${\mathbb C} \setminus \overline{\cal C}$.
With a slight abuse of notation, we call $\hat{\cal C}$ a \emph{cluster} and denote by $\hat{\cal C}_z$
the cluster containing $z$. We set $\hat{\cal C}_z = \emptyset$ if $z$ is not contained in any cluster $\hat{\cal C}$.

\begin{theorem} \label{thm:phase-trans}
Let ${\cal A}(\lambda,m,D)$ be a massive Brownian loop soup in $D$ with intensity $\lambda$ and mass function
$m$, and denote by ${\mathbb P}_{\lambda,m}$ the distribution of ${\cal A}(\lambda,m,{\mathbb C})$.
\begin{itemize}
\item If $\lambda>1/2$, $m$ is bounded and $D$ is bounded, with probability one the vacant set of
${\cal A}(\lambda,m,D)$ is totally disconnected.
\item If $\lambda \leq 1/2$ and $m$ is bounded away from zero, the vacant set of
${\cal A}(\lambda,m,{\mathbb C})$ contains a unique infinite connected component. Moreover, there is a
$\xi<\infty$ such that, for any $z \in {\mathbb C}$ and all $L>0$,
\begin{equation} \label{eq:exp-decay}
{\mathbb P}_{\lambda,m}(\text{\emph{diam}}(\hat{\cal C}_z) \geq L) \leq e^{-L/\xi} \; .
\end{equation}
\end{itemize}
\end{theorem}

Note that, although in a massive loop soup individual large loops are exponentially suppressed, 
Eq.~(\ref{eq:exp-decay}) is far from obvious, and in fact false when $\lambda>1/2$, since in that
equation the exponential decay refers to clusters of loops.
Theorem~\ref{thm:phase-trans} corresponds to Theorem 2.11 of \cite{camia-notes};
we omit the proof since it is not essential for this paper and it can be found in \cite{camia-notes}.

We note that there is some confusion in most of the existing literature 
regarding the critical intensity corresponding to the connectivity phase transition in the Brownian loop soup.
This is related to the choice of normalization of the Brownian loop measure and is connected to the discussion
in the footnote on p.~3.

It is shown in~\cite{sw} that, in the \emph{critical} case ($m=0$), if $D$ is bounded, the set of outer boundaries
of the clusters $\hat{\cal C}$ that are not surrounded by other outer boundaries are distributed like a Conformal
Loop Ensemble (CLE) in $D$. Hence, in the massive case, the set of outer boundaries of the clusters $\hat{\cal C}$
that are not surrounded by other outer boundaries can be thought of as a massive CLE in $D$.

\section{Random Walk Loop Soups} \label{sec:rwls}

The random walk loop soup was introduced by Lawler and Trujillo Ferreras \cite{ltf} as a lattice version
of the Brownian loop soup of Lawler and Werner \cite{lw}. More general versions, including ones with a
killing measure, have been extensively studied by Le Jan \cite{lejan1,lejan2} (see also~\cite{sznitman-notes}).

Let $k_x \geq 0$ for every $x \in {\mathbb Z}^2$ and define $p_{x,y} = 1/(k_x + 4)$ if $|x-y|=1$ and
$p_{x,y}=0$ otherwise. If $k_x=0$ for all $x$, $\{p_{x,y}\}_{y \in {\mathbb Z}^2}$ is the collection of
transition probabilities for the simple symmetric random walk on ${\mathbb Z}^2$. If $k_x \neq 0$, then
$p_{x,y} = \frac{1}{4}(1+\frac{k_x}{4})^{-1}<\frac{1}{4}$ and one can interpret
$\{p_{x,y}\}_{y \in {\mathbb Z}^2}$ as the collection of transition probabilities for a random walker
``killed'' at $x$ with probability $1-(1+\frac{k_x}{4})^{-1}=\frac{k_x}{k_x+4}$. (Equivalently, one can
introduce a ``cemetery'' state $\Delta$ not in ${\mathbb Z}^2$ to which the random walker jumps from
$x \in {\mathbb Z}^2$ with probability $\frac{k_x}{k_x+4}$, and where it stays forever once it reaches it.)
Because of this interpretation, we will refer to the collection ${\bf k}=\{k_x\}_{x\in{\mathbb Z}^2}$ as
\emph{killing rates}.


Given a $(2n+1)$-tuple $(x_0,x_1,\ldots,x_{2n})$ with $x_0=x_{2n}$ and $|x_i-x_{i-1}|=1$ for $i=1,\ldots,2n$,
we call \emph{rooted lattice loop} the continuous path $\tilde\gamma: [0,2n] \to {\mathbb C}$ with $\tilde\gamma(i)=x_i$
for integer $i=0,\ldots,2n$ and $\tilde\gamma(t)$ obtained by linear interpolation for other $t$. We call $x_0$ the
\emph{root} of the loop and denote by $|\tilde\gamma|=2n$ the \emph{length} or \emph{duration} of the loop.

Now let $D$ denote either ${\mathbb C}$ or a connected subset of ${\mathbb C}$. 
Following Lawler and Trujillo Ferreras \cite{ltf}, but within the more general framework of the previous paragraph,
we introduce the \emph{rooted random walk loop measure} $\nu^{r,{\bf k}}_D$ which assigns the loop $\tilde\gamma$ of
length $|\tilde\gamma|$, with root $x$, weight $|\tilde\gamma|^{-1} p_{x,x_1} p_{x_1,x_2} \ldots p_{x_{|\tilde\gamma|-1},x}$
if $x,\ldots,x_{|\tilde\gamma|-1} \in D$ and 0 otherwise.

The \emph{unrooted random walk loop measure} $\nu^{u,{\bf k}}_D$ is obtained from the rooted one by
``forgetting the root.'' More precisely, if $\tilde\gamma$ is a rooted lattice loop and $j$ a positive integer,
$\theta_j \tilde\gamma : t \mapsto \tilde\gamma(j+t \mod |\tilde\gamma|)$ is again a rooted loop.
This defines an equivalence relation between rooted loops; an \emph{unrooted lattice loop} is an equivalence
class of rooted lattice loops under that relation. By a slight abuse of notation, in the rest of the paper we will
use $\tilde\gamma$ to denote unrooted lattice loops and $\tilde\gamma(\cdot)$ to denote any rooted
lattice loop in the equivalence class of $\tilde\gamma$.
The $\nu^{u,{\bf k}}_D$-measure of the unrooted loop $\tilde\gamma$ is the sum of the $\nu^{r,{\bf k}}_D$-measures
of the rooted loops in the equivalence class of $\tilde\gamma$. 
The \emph{length} or \emph{duration}, $|\tilde\gamma|$, of an unrooted loop $\tilde\gamma$ is the length
of any one of the rooted loops in the equivalence class $\tilde\gamma$.


\begin{definition} \label{def:rwls}
A \emph{random walk loop soup} in $D$ with intensity $\lambda$ is a Poissonian realization from
$\lambda\nu^{u,{\bf k}}_D$.
\end{definition}

A realization of the random walk loop soup in $D$ is a multiset (i.e., a set whose elements can occur multiple times)
of unrooted loops. If we denote by $N_{\tilde\gamma}$ the multiplicity of $\tilde\gamma$ in a loop soup with
intensity $\lambda$, $\{N_{\tilde\gamma}\}$ is a collection of independent Poisson random variables with
parameters $\lambda\nu^{u,{\bf k}}_D(\tilde\gamma)$. Therefore, the probability that a realization of the
random walk loop soup in $D$ with intensity $\lambda$ contains each loop $\tilde\gamma$ in $D$ with
multiplicity $n_{\tilde\gamma}\geq0$ is equal to
\begin{equation} \label{rwls}
\prod_{\tilde\gamma} \exp{\left(-\lambda\nu^{u,{\bf k}}_D(\tilde\gamma)\right)}
\frac{1}{n_{\tilde\gamma}!} \left(\lambda \nu^{u,{\bf k}}_D(\tilde\gamma)\right)^{n_{\tilde\gamma}}
= \frac{1}{{\cal Z}_{\lambda,{\bf k}}} \prod_{\tilde\gamma} \frac{1}{n_{\tilde\gamma}!} \left(\lambda \nu^{u,{\bf k}}_D(\tilde\gamma)\right)^{n_{\tilde\gamma}} \, ,
\end{equation}
where the product $\prod_{\tilde\gamma}$ is over all unrooted lattice loops in $D$ and
\begin{equation} \label{partition-function}
{\cal Z}_{\lambda,{\bf k}} := \exp{\left(\lambda \sum_{\tilde\gamma} \nu^{u,{\bf k}}_D(\tilde\gamma)\right)}
= 1 + \sum_{n=1}^{\infty} \frac{1}{n!} \sum_{(\tilde\gamma_1,\ldots,\tilde\gamma_n)} \prod_{i=1}^n \lambda \nu^{u,{\bf k}}_D(\tilde\gamma_i) \, ,
\end{equation}
where the sum over $(\tilde\gamma_1,\ldots,\tilde\gamma_n)$ is over all ordered configurations of $n$ loops,
not necessarily distinct.
From a statistical mechanical viewpoint, ${\cal Z}_{\lambda,{\bf k}}$ can be interpreted as the grand canonical
partition function of a ``gas'' of loops, and one can think of the random walk loop soup as describing a grand
canonical ensemble of noninteracting loops (an ``ideal gas'') with the killing rates $\{k_x\}$ and the intensity
$\lambda$ as free ``parameters.''
(For more on the statistical mechanical interpretation of the model, see Sec.~6.4 of~\cite{bb}.)

When $k_x=0 \; \forall x \in D \cap {\mathbb Z}^2$, we use $\nu^u_D$ to denote the unrooted random walk
loop measure in $D$; for reasons that will be clear when we talk about scaling limits, later in this section,
a random walk loop soup obtained using such a measure will be called \emph{critical}.

Now let $m: {\mathbb C} \to {\mathbb R}$ be a nonnegative function; we say that a random walk loop soup
has \emph{mass (function)} $m$ if $k_x = 4(e^{m^2(x)}-1)$ for all $x \in D \cap {\mathbb Z}^2$, and call
\emph{massive} a random walk loop soup with mass $m$ that is not identically zero on $D \cap {\mathbb Z}^2$.
For a massive random walk loop soup in $D$ with intensity $\lambda$ and mass $m$ we use the notation
$\tilde{\cal A}(\lambda,m,D)$.

The next proposition gives a construction for generating a massive random walk loop soup from a critical one,
establishing a useful probabilistic coupling between the two (i.e., a way to construct the two loop soups on the
same probability space).
\begin{proposition} \label{coupling2}
A random walk loop soup in $D$ with intensity $\lambda$ and mass function $m$ can be realized in the following way.
\begin{enumerate}
\item Take a realization of the \emph{critical} random walk loop soup in $D$ with intensity $\lambda$.
\item Assign to each loop $\tilde\gamma$ an independent, mean-one, exponential random variable $T_{\tilde\gamma}$.
\item Remove from the soup the loop $\tilde\gamma$ of length $|\tilde\gamma|$ if
\begin{equation} \label{eq:loop-killing-2}
\sum_{i=0}^{|\tilde\gamma|-1} m^2(\tilde\gamma(i)) > T_{\tilde\gamma} \, .
\end{equation}
\end{enumerate}
\end{proposition}

\begin{remark}
Note that Eq.~\eqref{eq:loop-killing-2} requires choosing a rooted loop from the equivalence class
$\tilde\gamma$ but is independent of the choice.
\end{remark}

\noindent{\bf Proof.}
It is easy to see that the construction of Proposition~\ref{coupling2} defines a Poisson point process,
so we just need to compare the expected number of loops in $D$ generated by the construction
of Prop.~\ref{coupling2} with that of the massive random walk loop soup, to verify the relation
between killing rates $\{k_x\}$ and mass function $m$.

Writing $p_{x,y} = \frac{1}{k_x+4} = \frac{1}{4} \frac{4}{k_x+4}$ when $|x-y|=1$,
for the massive random walk soup we have
\begin{equation} \label{mean1}
\lambda \sum_{\tilde\gamma} \nu^{u,{\bf k}}_D(\tilde\gamma) =
\lambda \sum_{\tilde\gamma} \nu^u_D(\tilde\gamma) \prod_{i=0}^{|\tilde\gamma|-1} \frac{4}{k_{x_i}+4} \, ,
\end{equation}
where the sum $\sum_{\tilde\gamma}$ is over all unrooted loops in $D$ and the right hand side implies that we have
chosen a representative for $\tilde\gamma$ such that $\tilde\gamma(i)=x_i$, but is independent of the choice.

The expected number of loops resulting from the construction of Prop.~\ref{coupling2} is
\begin{eqnarray} \label{mean2}
\lefteqn{\lambda \sum_{\tilde\gamma} \nu^{u}_D(\tilde\gamma) \int_0^{\infty} e^{-t} 
\mathbbm{1}_{\{\sum_{i=0}^{|\tilde\gamma|-1} m^2(\tilde\gamma(i))<t\}} dt } \nonumber \\
& & = \lambda \sum_{\tilde\gamma} \nu^{u}_D(\tilde\gamma) \prod_{i=0}^{|\tilde\gamma|-1}e^{-m^2(x_i)} \, ,
\end{eqnarray}
where we have chosen the same representative with $\tilde\gamma(i)=x_i$ as before and the result is again
independent of the choice. Comparing Eqs.~\eqref{mean1} and~\eqref{mean2} gives
$e^{-m^2(x)}=\frac{4}{k_x+4}=(1+k_x/4)^{-1}$ or $k_x=4(e^{m^2(x)}-1)$. \fbox{}

\subsection{The Near-Critical Scaling Limit} \label{sec:scal-lim}

We are now going to consider scaling limits for the random walk loop soups defined above. 
Consider first a critical, full-plane, random walk loop soup
$\tilde{\cal A}_{\lambda} \equiv \tilde{\cal A}(\lambda,0,{\mathbb Z}^2)$.
Following~\cite{ltf}, for each integer $N \geq 2$, we define the \emph{rescaled loop soup}
\begin{equation} \label{rescaled-soup}
\tilde{\cal A}^N_{\lambda} = \{ \tilde\Phi_N \tilde\gamma : \tilde\gamma \in \tilde{\cal A}_{\lambda}\}
\text{ with } \tilde\Phi_N \tilde\gamma(t) = N^{-1} \tilde\gamma(2N^2t) \, .
\end{equation}
$\tilde\Phi_N \tilde\gamma$ is a lattice loop of duration $t_{\tilde\gamma}:=|\tilde\gamma|/(2N^2)$
on the rescaled lattice $\frac{1}{N}{\mathbb Z}^2$ and so
$\tilde{\cal A}^N_{\lambda}$ is a 
random walk loop soup on $\frac{1}{N}{\mathbb Z}^2$, with rescaled time.
It is shown in~\cite{ltf} that, as $N \to \infty$, $\tilde{\cal A}^N_{\lambda}$
converges to the Brownian loop soup of~\cite{lw} in an appropriate sense.
This means that the critical random walk loop soup has a conformally invariant scaling
limit (the Brownian loop soup), which explains our use of the term \emph{critical}.

If we rescale in the same way a massive random walk loop soup with constant mass function $m>0$,
the resulting scaling limit is trivial, in the sense that it does not contain any loops larger than one point.
This is so because, under the random walk loop measure, only loops of duration of order at least $N^2$
have diameter of order at least $N$ with non-negligible probability as $N \to \infty$, and are therefore
``macroscopic'' in the scaling limit. It is then clear that, in order to obtain a nontrivial scaling limit, the
mass function needs to be rescaled while taking the scaling limit.

Suppose, for simplicity, that the mass function $m$ is constant, and let $m_N$ denote the rescaled
mass function. When $m_N$ tends to zero, $k_x \approx 4 m_N^2$ and one has the following dichotomy.
\begin{itemize}
\item If $\lim_{N \to \infty}N m_N = 0$, loops with a number of steps of the order of $N^2$ or smaller
are not affected by the killing in the scaling limit and one recovers the critical Brownian loop soup.
\item If $\lim_{N \to \infty}N m_N = \infty$, all loops with a number of steps of the order of $N^2$ or
more are removed from the soup in the scaling limit and no ``macroscopic'' loop (larger than one point)
is left.
\end{itemize}

In view of this observation, it is possible to obtain a \emph{near-critical} scaling limit, that is, a nontrivial
scaling limit that differs from the critical one, only if the mass function $m$ is rescaled by $O(1/N)$.
This leads us to considering the 
loop soup $\tilde{\cal A}^N_{\lambda,m}$ defined as a random walk loop soup on the rescaled lattice
$\frac{1}{N}{\mathbb Z}^2$ with mass function $m/(\sqrt{2}N)$ and rescaled time as in~\eqref{rescaled-soup}.
Such a soup can be obtained from $\tilde{\cal A}^N_{\lambda}$ using the construction in
Prop.~\ref{coupling2}, replacing $m^2(\tilde\gamma(i))$ with $\frac{1}{2N^2}m^2(\tilde\gamma(i)/N)$
in Eq.~(\ref{eq:loop-killing-2}).

\begin{theorem} \label{prop:massive-soups}
Let $m$ be a nonnegative function such that $m^2$ is Lipschitz continuous.
There exist two sequences $\{{\cal A}^N_{\lambda,m}\}_{N\geq2}$ and $\{\tilde{\cal A}^N_{\lambda,m}\}_{N\geq2}$
of loop soups defined on the same probability space and such that the following holds.
\begin{itemize}
\item For each $\lambda>0$, ${\cal A}^N_{\lambda,m}$ is a 
massive Brownian loop soup in $\mathbb C$ with intensity $\lambda$ and mass $m$; the realizations
of the loop soup are increasing in $\lambda$.
\item For each $\lambda>0$, $\tilde{\cal A}^N_{\lambda,m}$ is a 
massive random walk loop soup on $\frac{1}{N}{\mathbb Z}^2$ with intensity $\lambda$, mass
$m/(\sqrt{2}N)$ and time scaled as in~\eqref{rescaled-soup}; the realizations of the loop soup are
increasing in $\lambda$.
\item For every bounded $D \subset {\mathbb C}$, with probability going to one as $N \to \infty$,
loops from ${\cal A}^N_{\lambda,m}$ and $\tilde{\cal A}^N_{\lambda,m}$ that are contained in $D$ and
have duration at least $N^{-1/6}$ can be put in a one-to-one correspondence with the following property. If
$\gamma \in {\cal A}^N_{\lambda,m}$ and $\tilde\gamma \in \tilde{\cal A}^N_{\lambda,m}$ are paired
in that correspondence and $t_{\gamma}$ and $t_{\tilde\gamma}$ denote their respective durations, then
\begin{eqnarray*}
|t_{\gamma} - t_{\tilde\gamma}| \leq \frac{5}{8} N^{-2} \\
\sup_{0 \leq s \leq 1} | \gamma(s t_{\gamma}) - \tilde\gamma(s t_{\tilde\gamma})| \leq c_1 N^{-1} \log N
\end{eqnarray*}
for some constant $c_1$.
\end{itemize}

\end{theorem}

\noindent{\bf Proof.} 
Let ${\cal A}_{\lambda}$ be a critical Brownian loop soup in $\mathbb C$ with intensity $\lambda$
and $\tilde{\cal A}_{\lambda}$ a critical random walk loop soup on ${\mathbb Z}^2$ with intensity
$\lambda$, coupled as in Theorem~1.1 of~\cite{ltf}. Consider the scaled loop soups ${\cal A}^N_{\lambda}$
and $\tilde{\cal A}^N_{\lambda}$, where $\tilde{\cal A}^N_{\lambda}$ is defined in~\eqref{rescaled-soup}
and ${\cal A}^N_{\lambda} := \{\Phi_N \gamma : \gamma \in {\cal A}_{\lambda}\}$ with
$\Phi_N \gamma(t) = N^{-1} \gamma(N^2t)$ for $0 \leq t \leq t_{\gamma}/N^2$.
Note that, because of scale invariance, ${\cal A}^N_{\lambda}$ is a critical Brownian loop soup
in $\mathbb C$ with parameter $\lambda$.

It readily follows from Theorem~1.1 of~\cite{ltf} that, if one considers only loops of duration
greater than $N^{-1/6}$, loops from ${\cal A}^N_{\lambda}$ and $\tilde{\cal A}^N_{\lambda}$ contained
in $D$ can be put in a one-to-one correspondence with the properties described in Theorem~\ref{prop:massive-soups},
except perhaps on an event of probability going to zero as $N \to \infty$. For simplicity, in the rest of the proof we will
call {\it macroscopic} the loops of duration greater than $N^{-1/6}$.

On the event that such a one-to-one correspondence between macroscopic loops in $D$ exists, we
construct the massive loop soups ${\cal A}^N_{\lambda,m}$ and $\tilde{\cal A}^N_{\lambda,m}$ in the
following way. To each pair of macroscopic loops $\gamma \in {\cal A}^N_{\lambda}$ and
$\tilde\gamma \in \tilde{\cal A}^N_{\lambda}$, paired in the correspondence of Theorem~1.1 of~\cite{ltf},
we assign an independent, mean-one, exponential random variable $T_{\gamma}$. We let $t_{\gamma}$
denote the (rescaled) duration of $\gamma$ and $t_{\tilde\gamma}$ the (rescaled) duration of $\tilde\gamma$,
and let $M=2N^2 t_{\tilde\gamma}$ denote the number of steps of the lattice loop $\tilde\gamma$.
As in the constructions described in Props.~\ref{coupling1} and \ref{coupling2}, we remove $\gamma$
from the Brownian loop soup if $\int_0^{t_{\gamma}} m^2(\gamma(s)) ds > T_{\gamma}$ and remove
$\tilde\gamma$ from the random walk loop soup if
$\frac{1}{2N^2}\sum_{k=0}^{M-1} m^2(\tilde\gamma(\frac{k}{2N^2})) > T_{\gamma}$.
The resulting loop soups, ${\cal A}^N_{\lambda,m}$ and $\tilde{\cal A}^N_{\lambda,m}$, are defined on the same
probability space and are distributed like a massive Brownian loop soup with mass function $m$ and a random
walk loop soup with mass function $m/(\sqrt{2}N)$, respectively. We use $\mathbb P$ to denote the joint distribution
of ${\cal A}^N_{\lambda,m}$, $\tilde{\cal A}^N_{\lambda,m}$ and the collection $\{T_{\gamma}\}$.

For loops that are not macroscopic, the removal of loops is done independently for the Brownian loop soup and the
random walk loop soup. If there is no one-to-one correspondence between macroscopic loops in $D$, the removal
is done independently for all loops, including the macroscopic ones.

We want to show that, on the event that there is a one-to-one correspondence between macroscopic loops in $D$,
the one-to-one correspondence survives the removal of loops described above with probability going to one as
$N \to \infty$. For that purpose, we need to compare $\int_0^{t_{\gamma}} m^2(\gamma(s)) ds$ and
$\frac{1}{2N^2}\sum_{k=0}^{M-1} m^2(\tilde\gamma(\frac{k}{2N^2}))$ for loops $\gamma$ and $\tilde\gamma$
paired in the above correspondence. In order to do that, we write
\begin{eqnarray*}
\int_0^{t_{\gamma}} m^2(\gamma(s)) ds & = & t_{\gamma} \int_0^1 m^2(\gamma(t_{\gamma}u)) du \\
& = & \lim_{n \to \infty} \frac{t_{\gamma}}{n t_{\tilde\gamma}} \sum_{i=0}^{\lfloor n t_{\tilde\gamma}\rfloor} m^2\left(\gamma\left(\frac{t_{\gamma} i}{n t_{\tilde\gamma}} \right)\right) \\
& = & \lim_{q \to \infty} \frac{t_{\gamma}/t_{\tilde\gamma}}{4qN^2} \sum_{i=0}^{2qM-1} m^2\left(\gamma\left(\frac{i}{2qM} t_{\gamma} \right)\right) \, ,
\end{eqnarray*}
where $t_{\tilde\gamma} = \frac{M}{2N^2}$ and the last expression is obtained by letting $n=4qN^2$, with $q \in {\mathbb N}$.
Thus, for fixed $N$ and $\gamma$, the quantity
\begin{equation} \nonumber
\Omega(N,q;\gamma) := \left| \int_0^{t_{\gamma}} m^2(\gamma(s)) ds - \frac{t_{\gamma}/t_{\tilde\gamma}}{4qN^2}
\sum_{i=0}^{2qM-1} m^2 \left(\gamma\left(\frac{i}{2qM} t_{\gamma} \right)\right) \right|
\end{equation}
can be made arbitrarily small by chossing
$q$ sufficiently large.

Define the sets of indeces $I_0 = \{ i: 0 \leq i < q \} \cup \{ i: (2M-1)q \leq i < 2qM \}$ and
$I_k = \{ i: (2k-1)q \leq i < (2k+1)q \}$ for $1 \leq k \leq M-1$.
For $i \in I_k$, $0 \leq k \leq M-1$, we have that
\begin{eqnarray*}
\lefteqn{\left|\gamma\left(\frac{i}{2qM} t_{\gamma}\right) - \tilde\gamma\left(\frac{k}{M}t_{\tilde\gamma}\right)\right| } \\
& & \leq \left|\gamma\left(\frac{i}{2qM} t_{\gamma}\right) - \tilde\gamma\left(\frac{i}{2qM}t_{\tilde\gamma}\right)\right|
+ \left|\tilde\gamma\left(\frac{i}{2qM}t_{\tilde\gamma}\right) - \tilde\gamma\left(\frac{k}{M}t_{\tilde\gamma}\right)\right| \\
& & \leq \frac{c_1 \log N}{N} + \frac{\sqrt 2}{N}
\end{eqnarray*}
for some constant $c_1$, where the first term in the last line comes from Theorem~1.1 of~\cite{ltf} and the second term
comes from the fact that $\tilde\gamma(s)$ is defined by interpolation and that
\begin{itemize}
\item if $i \in I_0$, either $0 \leq \frac{i}{2qM}t_{\tilde\gamma} < \frac{1}{2M}t_{\tilde\gamma}$ so that $\tilde\gamma(\frac{i}{2qM}t_{\tilde\gamma})$
falls on the edge of $\frac{1}{N} {\mathbb Z}^2$ between $\tilde\gamma(0)$ and $\tilde\gamma(\frac{t_{\tilde\gamma}}{M}) = \tilde\gamma(\frac{1}{2N^2})$,
or $(1 - \frac{1}{2M})t_{\tilde\gamma} \leq \frac{i}{2qM}t_{\tilde\gamma} < t_{\tilde\gamma}$
so that $\tilde\gamma(\frac{i}{2qM}t_{\tilde\gamma})$ falls on the edge between
$\tilde\gamma(t_{\tilde\gamma} - \frac{t_{\tilde\gamma}}{M}) = \tilde\gamma(t_{\tilde\gamma}-\frac{1}{2N^2})$ and $\tilde\gamma(t_{\tilde\gamma}) = \tilde\gamma(0)$,
\item if $i \in I_k$ with $1 \leq k \leq M-1$,
$\frac{k}{M}t_{\tilde\gamma} - \frac{t_{\tilde\gamma}}{2M} \leq \frac{i}{2qM}t_{\tilde\gamma} < \frac{k}{M}t_{\tilde\gamma} + \frac{t_{\tilde\gamma}}{2M}$
so that $\tilde\gamma(\frac{i}{2qM}t_{\tilde\gamma})$ falls either on the edge of $\frac{1}{N} {\mathbb Z}^2$ between
$\tilde\gamma(\frac{k-1}{M}t_{\tilde\gamma}) = \tilde\gamma(\frac{k-1}{2N^2})$ and $\tilde\gamma(\frac{k}{2N^2})$,
or on the edge between $\tilde\gamma(\frac{k}{2N^2})$ and $\tilde\gamma(\frac{k+1}{M}t_{\tilde\gamma}) = \tilde\gamma(\frac{k+1}{2N^2})$.
\end{itemize}

Since $m^2$ is Lipschitz continuous, for each $i \in I_k, \, 0 \leq k \leq M-1$, we have that
\begin{equation} \nonumber
\left|m^2\left(\gamma\left(\frac{i}{2qM}t_{\gamma}\right)\right) - m^2\left(\tilde\gamma\left(\frac{k}{M}t_{\tilde\gamma}\right)\right)\right| \leq \frac{c_2\log N}{N}
\end{equation}
for some constant $c_2<\infty$ and all $N\geq2$.
We let $\overline{m}^2_D:=\sup_{x \in D}m^2(x)$ and observe that, since $t_{\tilde\gamma} \geq N^{-1/6}$, the inequality
$|t_{\gamma}-t_{\tilde\gamma}|<\frac{5}{8}N^{-2}$ from Theorem~1 of~\cite{ltf} implies that
$t_{\gamma}/t_{\tilde\gamma} < 1 + \frac{5}{8} N^{-11/6}$ and that $M=2N^2 t_{\tilde\gamma} < 2N^2 t_{\gamma} + \frac{5}{4}$.
It then follows that
\begin{eqnarray*} 
\lefteqn{\left| \frac{1}{2N^2} \sum_{k=0}^{M-1} m^2\left(\tilde\gamma\left(\frac{k}{2N^2}\right)\right)
- \frac{t_{\gamma}/t_{\tilde\gamma}}{4qN^2} \sum_{i=0}^{2qM-1} m^2\left(\gamma\left(\frac{i}{2qM} t_{\gamma} \right)\right) \right| } \nonumber \\
& & \leq \frac{1}{2N^2} \sum_{k=0}^{M-1} \left| m^2\left(\tilde\gamma\left(\frac{k}{2N^2}\right)\right)
- \frac{t_{\gamma}/t_{\tilde\gamma}}{2q} \sum_{i \in I_k} m^2\left(\gamma\left(\frac{i}{2qM} t_{\gamma} \right)\right) \right| \nonumber \\
& & \leq \frac{1}{2N^2} \sum_{k=0}^{M-1} \frac{1}{2q} \sum_{i \in I_k} \left| m^2\left(\tilde\gamma\left(\frac{k}{2N^2}\right)\right)
- \frac{t_{\gamma}}{t_{\tilde\gamma}} m^2\left(\gamma\left(\frac{i}{2qM} t_{\gamma} \right)\right) \right| \nonumber \\
& & \leq  \frac{|1-t_{\gamma}/t_{\tilde\gamma}|}{2N^2} \sum_{k=0}^{M-1} m^2\left(\tilde\gamma\left(\frac{k}{2N^2}\right)\right) \nonumber \\
& & + \frac{t_{\gamma}/t_{\tilde\gamma}}{4qN^2} \sum_{k=0}^{M-1} \sum_{i \in I_k} \left| m^2\left(\tilde\gamma\left(\frac{k}{2N^2}\right)\right)
- m^2\left(\gamma\left(\frac{i}{2qM} t_{\gamma}\right)\right) \right| \nonumber \\
& & \leq  \frac{|1-t_{\gamma}/t_{\tilde\gamma}|}{2N^2} \sum_{k=0}^{M-1} m^2\left(\tilde\gamma\left(\frac{k}{2N^2}\right)\right) + t_{\gamma} \frac{c_2\log N}{N} \nonumber \\
& & \leq  \frac{|1-t_{\gamma}/t_{\tilde\gamma}|}{2N^2} M \overline{m}^2_D + t_{\gamma} \frac{c_2\log N}{N}
< t_{\gamma} \frac{c'_3\log N}{N} \, ,
\end{eqnarray*}
for some positive constant $c'_3=c'_3(D,m)<\infty$ independent of $\gamma$ and $\tilde\gamma$.

Therefore, for fixed $N$ and pair of macroscopic loops, $\gamma$ and $\tilde\gamma$, and for any $q \in {\mathbb N}$,
\begin{equation} \nonumber
\left| \int_0^{t_{\gamma}} m^2(\gamma(s)) ds -
\frac{1}{2N^2}\sum_{k=0}^{M-1} m^2\left(\tilde\gamma\left(\frac{k}{2N^2}\right)\right) \right|
< t_{\gamma} \frac{c'_3\log N}{N} + \Omega(N,q;\gamma) \, .
\end{equation}
For fixed $N$ and $\gamma$, one can choose $q^*$ so large that
\begin{equation} \nonumber
\Omega(N,q^*;\gamma) < t_{\gamma} \frac{c'_3 \log N}{N} \, .
\end{equation}
Hence, there is a positive constant $c_3=2c'_3$ such that, for every $N \geq 2$ and every pair of
macroscopic loops, $\gamma$ and $\tilde\gamma$, paired in the correspondence of Theorem~1.1 of~\cite{ltf},
\begin{equation} \nonumber
\left| \int_0^{t_{\gamma}} m^2(\gamma(s)) ds -
\frac{1}{2N^2}\sum_{k=0}^{M-1} m^2\left(\tilde\gamma\left(\frac{k}{2N^2}\right)\right) \right|
< t_{\gamma} \frac{c_3\log N}{N} \, .
\end{equation}

We now need to estimate the number of macroscopic loops contained in $D$. For that purpose,
we note that, using the rooted Brownian loop measure~\eqref{eq:rooted-brownian-loop-measure},
the mean number, $\cal M$, of macroscopic loops contained in $D$ can be bounded above by
\begin{equation} \label{eq:mean}
{\cal M} = \lambda \int_D \int_{N^{-1/6}}^{\infty} \frac{1}{2\pi t^2} \, \mu^{br}_{z,t}(\gamma: \gamma \subset D) \, dt \, dA(z) \leq \frac{\lambda \, \diam^2(D)}{8} N^{1/6} \, .
\end{equation}

Let ${\cal A}^N(\lambda,m;D)$ (respectively, $\tilde{\cal A}^N(\lambda,m;D)$) denote the massive Brownian
(resp., random walk) loop soup in $D$, i.e., the set of loops from ${\cal A}^N_{\lambda,m}$
(respectively, $\tilde{\cal A}^N_{\lambda,m}$) contained in $D$. For the critical soups, we use the
same notation omitting the $m$.

Let $A_N$ denote the event that there is a one-to-one correspondence between macroscopic loops from
${\cal A}^N(\lambda;D)$ and $\tilde{\cal A}^N(\lambda;D)$, and let $A^m_N$ denote the event that
there is a one-to-one correspondence between macroscopic loops from ${\cal A}^N(\lambda,m;D)$ and
$\tilde{\cal A}^N(\lambda,m;D)$. Furthermore, we denote by $X$ the number of macroscopic loops in
${\cal A}^N(\lambda;D)$, and by $T$ a mean-one exponential random variable. We have that, for any
$c_4,\theta>0$ and for all $N$ sufficiently large,

\begin{eqnarray*}
{\mathbb P}(A^m_N) & \geq & {\mathbb P}(A^m_N \cap A_N \cap \{X \leq c_4 N^{1/6}\}
\cap \{ \nexists \gamma \in {\cal A}^N(\lambda;D) : t_{\gamma} \geq \theta \}) \\
& = & {\mathbb P}(A^m_N | A_N \cap \{X \leq c_4 N^{1/6}\}
\cap \{ \nexists \gamma \in {\cal A}^N(\lambda;D) : t_{\gamma} \geq \theta \}) \\
& & {\mathbb P}(A_N \cap \{X \leq c_4 N^{1/6}\}
\cap \{ \nexists \gamma \in {\cal A}^N(\lambda;D) : t_{\gamma} \geq \theta \}) \\
& \geq & \left[ 1 - \sup_{x \geq 0} \Pr\left(x \leq T \leq x + \frac{c_3 \theta \log N}{N}\right)
\right]^{c_4 N^{1/6}} \\
& & {\mathbb P}(A_N \cap \{X \leq c_4 N^{1/6}\}
\cap \{ \nexists \gamma \in {\cal A}^N(\lambda;D) : t_{\gamma} \geq \theta \}) \\
& = & \exp\left(-\frac{c_5 \theta \log N}{N^{5/6}}\right) \\
& & {\mathbb P}(A_N \cap \{X \leq c_4 N^{1/6}\}
\cap \{ \nexists \gamma \in {\cal A}^N(\lambda;D) : t_{\gamma} \geq \theta \}) \, ,
\end{eqnarray*}
where $c_5=c_3 c_4$.

Since $\exp\left(-\frac{c_5 \theta \log N}{N^{5/6}}\right) \to 1$ as $N \to \infty$ for any $c_5,\theta>0$,
in order to conclude the proof, it suffices to show that ${\mathbb P}(A_N \cap \{ X \leq c_4 N^{1/6} \} \cap
\{ \nexists \gamma \in {\cal A}^N(\lambda;D) : t_{\gamma} \geq \theta \})$ can be made arbitrarily close
to one for some choice of $c_4$ and $\theta$, and $N$ sufficiently large.
But by Theorem~1 of~\cite{ltf}, ${\mathbb P}(A_N) \geq 1 - c (\lambda+1) \diam^2(D) N^{-7/2} \to 1$ as $N \to \infty$;
moreover, if $c_4>\frac{\lambda \, \diam^2(D)}{8}$, by Eq.~\eqref{eq:mean}, $c_4 N^{1/6}$ is larger than
the mean number of macroscopic loops in $D$. Since $X$ is a Poisson random variable with parameter equal to
the mean number $\cal M$ of macroscopic loops in $D$, the latter fact (together with a Chernoff bound argument)
implies that
\begin{equation} \nonumber
{\mathbb P}(X > c_4 N^{1/6}) \leq \frac{e^{-{\cal M}} (e {\cal M})^{c_4 N^{1/6}}}{(c_4 N^{1/6})^{c_4 N^{1/6}}}
\leq \left(\frac{e \, \lambda \, \diam^2(D)}{8 c_4}\right)^{c_4 N^{1/6}} \, .
\end{equation}
This shows that, if $c_4 > e \, \lambda \, \diam^2(D)/8$, ${\mathbb P}(X \leq c_4 N^{1/6}) \to 1$ as $N \to \infty$.

To find a lower bound for ${\mathbb P}(\nexists \gamma \in {\cal A}^N(\lambda;D) : t_{\gamma} \geq \theta)$, we define
\begin{equation} \nonumber
{\cal L}_{\theta,D} := \{ \text{loops } \gamma \text{ with } t_{\gamma} \geq \theta \text{ that stay in } D \} \, .
\end{equation}
We then have
\begin{eqnarray} \label{lower-bound}
{\mathbb P}(\nexists \gamma \in {\cal A}^N(\lambda;D) : t_{\gamma} \geq \theta) & = &
\exp\left[ -\lambda \mu_D({\cal L}_{\theta,D}) \right] \nonumber \\
& \geq & 1 - \lambda \mu_D({\cal L}_{\theta,D}) \nonumber \\
& \geq & 1 - \frac{\lambda \, \diam^2(D)}{\theta} \, ,
\end{eqnarray}
where the last line follows from the bound
\begin{equation} \nonumber
\mu_D({\cal L}_{\theta,D})
= \int_D \int_{\theta}^{\infty} \frac{1}{2 \pi t^2} \, \mu^{br}_{z,t}(\gamma : \gamma \text{ stays in } D) \, dt \, d{\bf A}(z)
\leq \frac{\diam^2(D)}{\theta} \, .
\end{equation}

The lower bound~\eqref{lower-bound}, together with the previous observations, shows that ${\mathbb P}(A_N^m)$
can be made arbitrarily close to one by choosing $c_4 > e \, \lambda\, \diam^2(D)/8$, $\theta$ sufficiently large,
depending on $D$, and then $N$ sufficiently large, depending on the values of $c_4$ and $\theta$. \fbox{}

\section{Boundary Perturbations--Proof of Thm.~\ref{theorem-rd-derivative}} \label{sec:boundary}

Let $D \subset {\mathbb C}$ be a bounded, simply connected domain with two marked points, $a,b \in \partial D$, on
its boundary. Let $D^{\#} := D \cap {\mathbb Z}^2$ and assign to it a graph structure by declaring two vertices,
$x,y \in D^{\#}$, \emph{adjacent} if $|x-y|=1$ and the straight segment between $x$ and $y$ is contained in $D$.
A vertex $x \in D^{\#}$ is on the \emph{boundary} of $D^{\#}$ if it is at distance $1$ from a vertex of ${\mathbb Z}^2$
that is not in $D^{\#}$ or from a vertex of $D^{\#}$ to which it is not adjacent. Choose two points $a^{\#}, b^{\#}$ on
the boundary of $D^{\#}$ so that they belong to the same connected component of $D^{\#}$ and are as close as possible
to $a$ and $b$, respectively.

Let $\tilde\omega_n=(x_0,x_1,\ldots,x_n)$ be a path such that $x_{i} \in D^{\#}$ and $x_{i}$ and $x_{i-1}$ are adjacent
for all $i=1,\ldots,n$, and $x_{0} = a^{\#}, x_{n} = b^{\#}$. Let $LE(\tilde\omega)$ denote the \emph{chronological loop
erasure} of $\tilde\omega$ (see, e.g., Section~9.5 of \cite{lawler-limic-book}), which generates a self-avoiding path
$\tilde\eta=LE(\tilde\omega)$ in $D^{\#}$ from $a^{\#}$ to $b^{\#}$. If $\tilde\omega$ is not self-avoiding, the
chronological loop erasure procedure also generates an ordered collection of \emph{rooted} loops
$\Gamma(\tilde\omega) = (\tilde\gamma^r_{i}(\tilde\omega))_{i=1,\ldots,k(\tilde\omega)}$.
(In this section, we use the notation $\tilde\gamma^r$ to distinguish rooted loops from unrooted ones.)

To a path $\tilde\omega$ we assign weight
\begin{equation} \nonumber
p_{x_{0},x_{1}} p_{x_{1},x_{2}} \ldots p_{x_{n-1},x_{n}} = \prod_{x \in \tilde\omega} (k_{x}+4)^{-n(x,\tilde\omega)},
\end{equation}
where $n(x,\tilde\omega)$ denotes the number of visits of $\tilde\omega$ to $x$.
We now assign to a self-avoiding path $\tilde\eta$ a weight $\nu^{\bf k}_{D}(\tilde\eta)$ defined as
\begin{equation*}
\nu^{\bf k}_{D}(\tilde\eta) :=
\sum_{\tilde\omega: LE(\tilde\omega)=\tilde\eta} \prod_{x \in \tilde\omega} (k_{x}+4)^{-n(x,\tilde\omega)}.
\end{equation*}
Splitting $\tilde\omega$ into $\tilde\eta=LE(\tilde\omega)$ and the ordered collection of rooted loops
$\Gamma(\tilde\omega) = (\tilde\gamma^r_{i}(\tilde\omega))_{i=1,\ldots,k(\tilde\omega)}$, and using Proposition~9.5.1
of~\cite{lawler-limic-book}, we obtain
\begin{eqnarray*}
\nu^{\bf k}_{D}(\tilde\eta) & = & \prod_{x \in \tilde\eta} (k_{x}+4)^{-n(x,\tilde\eta)} \, \Big( 1 +
\sum_{\tilde\omega \neq \tilde\eta: LE(\tilde\omega)=\tilde\eta} \prod_{i=1}^{k(\tilde\omega)}
\prod_{x \in \tilde\gamma^r_{i}(\tilde\omega)} (k_{x}+4)^{-n(x,\tilde\gamma^r_{i}(\tilde\omega))} \Big) \\
& = & \prod_{x \in \tilde\eta} (k_{x}+4)^{-n(x,\tilde\eta)} \, \exp\Big(
\sum_{\tilde\gamma^r \subset D: \tilde\gamma^r \cap \tilde\eta \neq \emptyset} \frac{1}{|\tilde\gamma^r|}
\prod_{x \in \tilde\gamma^r} (k_{x}+4)^{-n(x,\tilde\gamma^r)} \Big) ,
\end{eqnarray*}
where the last sum is over all possible rooted loops in $D$ intersecting $\tilde\eta$.

If $D' \subset D$ is a simply-connected subset of $D$ such that $\partial D'$ agrees with $\partial D$ near
$a$ and $b$, let $D_1 = D \setminus D'$. Then, for any self-avoiding path $\tilde\eta$ in $(D')^{\#}$ from
$a^{\#}$ to $b^{\#}$, we have that
\begin{eqnarray*}
\frac{\nu^{\bf k}_{D'}(\tilde\eta)}{\nu^{\bf k}_{D}(\tilde\eta)} & = & \exp\Big(
-\sum_{\tilde\gamma^r \subset D: \tilde\gamma^r \cap \tilde\eta \neq \emptyset, \tilde\gamma^r \cap D_1 \neq \emptyset}
\frac{1}{|\tilde\gamma^r|} \prod_{x \in \tilde\gamma^r} (k_{x}+4)^{-n(x,\tilde\gamma^r)} \Big) \\
& = & \exp\Big(-\nu_{D}^{u,\bf k}\left(\tilde\gamma \subset D: \tilde\gamma \cap \tilde\eta \neq \emptyset, \tilde\gamma \cap D_1 \neq \emptyset \right) \Big) \\
& = & \tilde{\mathbb P}^{D}_{1,{\bf k}}\left( \text{no loop intersects } D_1 \text{ and } \tilde\eta \right) ,
\end{eqnarray*}
where $\tilde{\mathbb P}^{D}_{\lambda,{\bf k}}$ denotes the distribution of a random walk loop soup in $D$ with
intensity $\lambda$ and killing rates $\bf k$.

To obtain a probability measure, $\tilde{\cal P}^{\bf k}_{D}$, associated to the loop-erased random walk with killing, one needs
to normalize $\nu^{\bf k}_{D}$ by its mass,
\begin{equation} \nonumber
\left| \nu^{\bf k}_{D} \right| = \sum_{\tilde\eta} \prod_{x \in \tilde\eta} (k_x+4)^{-n(x,\tilde\eta)} \,
e^{\nu^{u,\bf k}_{D}(\tilde\gamma: \tilde\gamma \cap \tilde\eta \neq \emptyset)} ,
\end{equation}
where the sum is over all self-avoiding paths $\tilde\eta$ in $D^{\#}$ from $a^{\#}$ to $b^{\#}$. Note that
\begin{equation} \nonumber
\frac{\left| \nu^{\bf k}_{D'} \right|}{\left| \nu^{\bf k}_{D} \right|} = \sum_{\tilde\eta \in D'}
\tilde{\cal P}^{\bf k}_{D}(\tilde\eta) \, e^{-\nu^{u,\bf k}_{D}(\tilde\gamma: \tilde\gamma \cap \tilde\eta \neq \emptyset,
\tilde\gamma \cap D_1 \neq \emptyset)} .
\end{equation}
Then, introducing
\begin{equation} \nonumber
\tilde Y_{\bf k} (\tilde\eta) = \tilde Y_{\bf k}(D,D';a,b)(\tilde\eta) :=
\tilde{\mathbb P}^{D}_{1,{\bf k}}(\text{no loop intersects } D_1 \text{ and } \tilde\eta) \mathbbm{1}_{\{ \tilde\eta \subset D' \}}
\end{equation}
and letting $\tilde E_{D}^{\bf k}$ denote expectation with respect to $\tilde{\cal P}_{D}^{\bf k}$, we have that
\begin{equation} \label{lattice-rd-derivative}
\frac{d\tilde{\cal P}_{D'}^{\bf k}}{d\tilde{\cal P}_{D}^{\bf k}}(\tilde\eta) = \frac{\tilde Y_{\bf k}(\tilde\eta)}{\tilde E_{D}^{\bf k}(\tilde Y_{\bf k})} .
\end{equation}

We now take $m \in {\mathbb R}$, $m \neq 0$, replace ${\mathbb Z}^2$ by $\frac{1}{N}{{\mathbb Z}^2}$, and take the scaling limit,
$N \to \infty$, of the random walk loop soup (in the sense of Theorem~\ref{prop:massive-soups}) and of the loop-erased random walk,
both with killing rates $k^{N}_x=4\left(e^{m^2/(2N^2)}-1 \right) \approx \frac{2m^2}{N^2}$ for every $x \in {\mathbb Z}^2$.
The random walk loop soup converges to the Brownian loop soup with mass $m$ by Theorem~\ref{prop:massive-soups}, while Theorem~2.1
of \cite{ms10} guarantees that the loop-erased random walk converges to the massive SLE$_2$ of \cite{ms10}. This observation, combined
with \eqref{lattice-rd-derivative}, concludes the proof of Eq.~\eqref{intro-rd-derivative} and of Theorem~\ref{theorem-rd-derivative}. \fbox{}\\


\noindent {\bf Acknowledgements.} The author thanks 
Tim van de Brug and Marcin Lis for comments on earlier drafts of the paper and for numerous discussions.
At the time of the research, the author was partially supported by NWO grant VIDI 639.032.916.

\end{document}